\documentclass[]{amsart}

\newtheorem{theorem}{Theorem}[section]
\newtheorem{lemma}[theorem]{Lemma}
\newtheorem{proposition}[theorem]{Proposition}

\newtheorem{corollary}[theorem]{Corollary}

\theoremstyle{definition}

\newtheorem{notation}[theorem]{Notation}

\def\owb{[\kern-0.2em[}
\def\cwb{]\kern-0.2em]}

\newcommand{\s}[2]{S\genfrac{[}{]}{0pt}{}{#1}{#2}}
\newcommand{\spr}[2]{S'\genfrac{[}{]}{0pt}{}{#1}{#2}}
\newcommand{\st}[2]{\widetilde{S}\genfrac{[}{]}{0pt}{}{#1}{#2}}

\newcommand{\ts}[2]{T\genfrac{[}{]}{0pt}{}{#1}{#2}}
\newcommand{\tpr}[2]{T'\genfrac{[}{]}{0pt}{}{#1}{#2}}

\let\emptyset=\mynull

\newcommand{\ab}{\allowbreak}

\title[Almost Mathieu Operator]
{On the Characteristic Polynomial \\ [5pt]
of the Almost Mathieu Operator}

\author[Lamoureux]{Michael P. Lamoureux}
\address{Department of Mathematics and Statistics, 
University of Calgary,
Calgary, Alberta, T2T 1A1}            

\email{mikel@math.ucalgary.ca}

\author[Mingo]{James A. Mingo}
\address{Department of Mathematics and
  Statistics, Queen's University,  King\-ston, Ontario, K7L 3N6}

\email{mingo@mast.queensu.ca}

\thanks{Research supported by a Discovery Grant 
from the Natural Sciences and Engineering
Research Council of Canada, AMS Classification: 47B 39 
(47B 15, 46L 05)}

\date{}

\begin{document}

\begin{abstract}
Let $A_\theta$ be the rotation C*-algebra for angle $\theta$. For
$\theta = p/q$ with $p$ and $q$ relatively prime, $A_\theta$ is the
sub-C*-algebra of $M_q(C(\mathbb{ T}^2))$ generated by a pair of
unitaries $u$ and $v$ satisfying $uv = e^{2 \pi i \theta} v u$.
Let $h_{\theta, \lambda} = u + u^{-1} + \lambda/2(v + v^{-1})$ be
the almost Mathieu operator. By proving an identity of rational
functions we show that for $q$ even, the constant term in the characteristic
polynomial of $h_{\theta, \lambda}$ is $(-1)^{q/2}(1 + (\lambda/2)^q) - (z_1^q
+ z_1^{-q} + (\lambda/2)^q(z_2^q + z_2^{-q}))$. 
\end{abstract}

\maketitle

\section{Introduction}

Let $\theta$, $\lambda$, and $\psi$ be
real numbers with $\lambda$ positive. The second order difference
operator $H_{\theta, \lambda,
\psi}$ on $\ell^2(\mathbb{ Z})$ given by
\[H_{\theta, \lambda, \psi}(\xi)(n) = \xi(n+1) + \xi(n-1) + \lambda
\cos(2 \pi n \theta + \psi)\xi(n) \]
for $\xi \in \ell^2(\mathbb{ Z})$ is called the almost Mathieu operator.
$H_{\theta, \lambda, \psi}$ is a discrete Schr\"odinger operator which
models an electron moving in a crystal lattice in a plane perpendicular
to a magnetic field.

An object of much study has been the spectrum $\sigma(\theta, \lambda) =
\cup_\psi \sigma(H_{\theta,\lambda,\psi})$. In \cite{hofstadter},
Hofstadter calculated $\sigma(\theta,2)$ for $\theta = p/q$ and $1 \leq p
< q\leq 50$. The remarkable pattern he found is called Hofstadter's
butterfly. For irrational $\theta$, a long standing concern has been the
connectedness and Lebesgue measure of $\sigma(\theta,\lambda)$  and the
labelling of the gaps, about which quite a bit is now known (see
\cite{avilaj}, \cite{avilak}, and \cite{puig} for spectacular recent 
advances as well as  \cite{avms}, \cite{bellissard}, \cite{boca},
\cite{choi-elliott-yui}, \cite{last} for earlier work). In addition there
has been numerical work on computing the spectrum to high accuracy for
large $q$
\cite[{\sc a$_2$, l}]{arveson1}.

Let $A_\theta$ be the rotation C*-algebra (see \cite{boca}). For $\theta =
p/q$ with $p$ and $q$ relatively prime and $\rho = e^{2 \pi \theta}$ let 
\[
 u_\theta = 
\left(\begin{matrix}
0  & 1 &       &        & 0  \\
   & 0 &   1   &        &   \\
   &   &   0   & \ddots  &   \\
   &   &       & \ddots & 1 \\
1  &   &       &        & 0 \\
\end{matrix}\right) \mbox{ and }
v_\theta = 
\left(\begin{matrix}
\rho  &        &        &            &   \\
      & \rho^2 &        &            &   \\
      &        & \ddots &            &   \\
      &        &        & \rho^{q-1} &   \\
      &        &        &            & 1 \\
\end{matrix}\right)
\]
i.e $u_\theta$ cyclically permutes the elements of the standard
basis and $v_\theta$ is a diagonal operator. Then define $u, v:
\mathbb{ T}^2 \rightarrow M_q(\mathbb{ C})$ by $u(z_1, z_2) = z_1
u_\theta$ and $v(z_1, z_2) = z_2 v_\theta$. Then $u\, v = \rho v\,
u$ and $A_\theta$ is the C*-algebra generated by $u$ and $v$ (see
\cite{boca}). The operator $h_{\theta, \lambda}
= u + u^{-1} + \lambda/2(v + v^{-1})$ contains all the spectral information of
$H_{\theta, \lambda, \psi}$ in that ${\rm Sp}(h_{\theta, \lambda}) =
\sigma_{\theta, \lambda} := \cup_\psi {\rm Sp}(H_{\theta, \lambda, \psi})$.

The main tool in the analysis of $\sigma_{\theta, \lambda}$ is
$\Delta_{\theta, \lambda}$, the discrete analogue of the discriminant. For
$\theta = p/q$, $\Delta_{\theta, \lambda}(x) = {\rm Tr}(A_1(x) \cdots
A_q(x))$ where 
$$ A_k(x) = 
\left(\begin{matrix}
 x - \lambda \cos(2 \pi k p/q + \pi/(2q)) & -1 \\
                      1                             &  0 \\
\end{matrix}\right)$$

Below are the first few values of this polynomial. Note the form of
$\Delta_{\theta, \lambda}$ so displayed depends only on the denominator
$q$; however, $\xi_\theta = 2 \cos(2 \pi p/q)$ depends on
the numerator $p$.

{\small
\renewcommand{\arraystretch}{1.5}
$\begin{array}{l|l}
q & \Delta_{\theta, 2}(x)
\mbox{ for } \theta= p/q \mbox{ and }
               \xi_\theta = 2 \cos(2 \pi \theta) 
\\ \hline 
2 & x^2 -4 \\
3 & x^3 -6  x \\
4 & x^4 - 8  x^2 + 4 \\
5 & x^5 - 10  x^3 + 5(3 - \xi_\theta)  x\\
6 & x^6 - 12  x^4 + 6(5 - \xi_\theta)  x^2 -4 \\
7 & x^7 -14  x^5 + 7(7 - \xi_\theta)   x^3 - 7(6 - 2\xi_\theta
    + 2\xi_{2\theta})  x \\
8 & x^8 - 16  x^6 + 8(9 - \xi_\theta)  x^4 - 8(12 - 4
\xi_\theta
    + 2 \xi_{2\theta})  x^2 + 4 \\
9 & x^9 - 18  x^7 + 9(11 - \xi_\theta)   x^5 - 9(31/3 -
6\xi_\theta
    + 2 \xi_{2 \theta})  x^3 + 9(14 - 8\xi_\theta + 3 \xi_{2
    \theta})  x \\
\end{array}$
}

One can calculate for $k=1, 2$ the coefficient of $x^{q - 2k}$, for $k =3$
the formula is conjectural (from numerical evidence). A deeper
understanding of the structure of $\Delta_{\theta. \lambda}$ would be
quite interesting. 

{\small\renewcommand{\arraystretch}{1.5}$
\begin{array}{l|l}
k & \mbox{coefficient \ of } x^{q - 2k} \mbox{ in } \Delta_{\theta,
\lambda}\quad (\mu = \lambda/2) \\ \hline
1 & -q(1 + \mu^q) \\
2 & q(\frac{1}{q-2}\binom{q-2}{2} \mu^4 + (q - 4 - \xi_\theta )
\mu^2 + \frac{1}{q-2}\binom{q-2}{2} \\ 
3 & -q( \frac{1}{q-3} \binom{q - 3}{3} \mu^6 +
(1 + \binom{q - 5}{2} - (q-6) \xi_\theta + \xi_{2 \theta}) \mu^4 \\
 & \mbox{} +
(1 + \binom{q - 5}{2} - (q-6) \xi_\theta + \xi_{2 \theta}) \mu^2 +
\frac{1}{q-3} \binom{q - 3}{3} )
 \end{array}$}

\

The connection with the characteristic polynomial of $h_{\theta, \lambda}$ is
given by
\begin{equation}\label{eqn00}
\det(xI_q - h_{\theta, \lambda}(z_1, z_2)) = \Delta_{\theta, \lambda}(x) +
z_1^q + z_1^{-q} + (\lambda/2)^q (z_2^q + z_2^{-q})
\end{equation}
and thus $\sigma_{\theta, \lambda} = \Delta_{\theta, \lambda}^{-1}[
-2(1 + (\lambda/2)^q),\ 2(1 + (\lambda/2)^q) ] $. Indeed,
$\Delta_{\theta, \lambda}(x)$ can be written as a determinant (c.f. Toda
\cite[\S 4]{toda})

\begin{equation}\label{eqn0}
\Delta_{p/q, \lambda}(x) =
\det\left(\begin{matrix}
\alpha_1  & 1   &        &        & 1   \\
1    & \alpha_2 & 1      &        &   \\
     & 1   & \ddots & \ddots &       \\
     &     & \ddots & \ddots &   \\
     &     &        &        & 1   \\
1    &     &        &      1 & \alpha_q \\
\end{matrix}\right) + 2 \bigg\{ (-1)^q + (\lambda/2)^q \bigg\}
\end{equation}
where all the other entries are 0 and  $\alpha_k = x -  \lambda \cos(2
\pi k p/q + \pi/(2q))$. Since  \begin{equation}\label{eqn1}
\Delta_{p/q, \lambda}(-x) = (-1)^q \Delta_{p/q, \lambda}(x)
\end{equation}
the coefficient of $x^{q - (2k+1)}$ is 0 for $0 \leq k < q/2$.

The main result of the paper asserts that for $a_l = 2 \cos(2 \pi l p/q)$ and
$1 \leq k < q/2$ we have
$$
\sum_{i_1, i_2, \dots , i_{q - 2k}} a_{i_1} a_{i_2} \cdots a_{i_{q - 2k}}
= 0$$ 
where the summation is over all subsets of $\{1, 2, 3, \dots , q\}$ obtained
by deleting $k$ pairs of adjacent elements -- counting 1 and $k$ as adjacent.
This is proved by establishing the following identity for $k \geq 3$ and $q
\geq 2k-1$
$$
\sum_{i_1 = 1}^{q - 2(k-1)} 
\kern-5pt\cdots \kern-7pt
\sum_{i_k = i_{k-1}+2}^q \ 
\prod_{j=1}^k \frac{(x^{-i_j} + x^{i_j})^{-1}}{%
              (x^{-i_j-1} + x^{i_j+1})  } 
=
\frac{\displaystyle
(x^{-q} - x^{q}) \prod_{i=k+1}^{2k-2} (x^{-q+i} - x^{q-i})}
{\displaystyle
\prod_{i=1}^k (x^{-2i} - x^{2i}) \prod_{i=-1}^{k-2}
(x^{-q+i} + x^{q-i}) }  
$$ 
$$\mbox{} +
\frac{(x^{-1} + x^1)^{-1} (x^{-q} + x^q)^{-1}}%
{(x^{-2} + x^2)  (x^{-q-1} + x^{q+1})}
\sum_{i_1 = 3}^{q - 2(k-2)} 
\kern-5pt\cdots \kern-5pt
\sum_{i_{k-2} = i_{k-3}+2}^{q-2} 
\prod_{j=1}^{k-2} \frac{ (x^{-i_j} + x^{i_j})^{-1}}{%
                    (x^{-i_j-1} + x^{i_j+1})       } 
$$

We then use this to show that for $a_l = 2 \cos(2 \pi l p /q)$ 
$$
\det\left(\begin{matrix}
a_1  & 1      &        &        & 1   \\
1    & a_2    & \ddots &        &   \\
     & \ddots & \ddots & \ddots &       \\
     &        & \ddots & \ddots & 1 \\
1    &        &        & 1      & a_q   \\
\end{matrix}\right)
= \begin{cases} 0  & q  \equiv 0   \pmod 4 \\
          4  & q  \equiv 1,3 \pmod 4 \\
          -8 & q  \equiv 2   \pmod 4 \\
\end{cases}
$$
  From this we show that the constant term (i.e. the coefficient of $x^0$) in
$\det(xI_q - h_{\theta, \lambda}(z_1, z_2))$ is
$$ (-1)^{q/2} 2 (1 + (\lambda/2)^q)) - (z_1^q + z_1^{-q} + 
(\lambda/2)^q (z_2^q + z_2^{-q}))$$
when $q$ is even. When $q$ is odd it
follows from (\ref{eqn1}) that the coefficient of $x^0$ is $ - (z_1^q +
z_1^{-q} + (\lambda/2)^q ( z_2^q + z_2^{-q}))$. 

Similar, though simpler, reasoning shows that the coefficient of $x^{q-2}$ is
$- q ( 1 + \lambda/2)$ and that the coefficient of $x^{q-4}$ is $
(\lambda/2)^4 q (q-3)/2 + (\lambda/2)^2 q (q - 4 - 2 \cos( 2 \pi \theta) ) +
q (q-3)/2$.

\section{The Main Theorem}

Let us use the following notation: let $a_1, \dots , a_n$ be elements of
a commutative ring and let
$$
(a_1, a_2, \dots , a_n) =
\left|
\begin{matrix}
a_1 & 1   &             &          & 0 \\
1   & a_2 & 1 &         &          &   \\
    &  \ddots &   a_3   & \ddots   &   \\
    &         & \ddots  & \ddots   & 1 \\
0   &         &         &  1 & a_n \cr \end{matrix}
\right|
$$
and
$$
\owb a_1, a_2, \dots , a_n\cwb =
\left|\begin{matrix}
a_1 & 1   &             &          & 1 \\
1   & a_2 & 1 &         &          &   \\
    &  \ddots &   a_3   & \ddots   &   \\
    &         & \ddots  & \ddots   & 1 \\
1   &         &         &  1 & a_n \cr \end{matrix}
\right|.
$$
The first matrix is a tridiagonal matrix with 1's on the sub and
super-diagonal and 0's elsewhere. The second matrix is the same tridiagonal
matrix with in addition 1's in the upper right and lower left corners, all
other entries are 0. Expanding along the bottom row we have

\begin{equation}\label{eqn2}
\owb a_1 , a_2, \dots , a_n\cwb = (a_1, a_2, \dots , a_n) 
- (a_2, a_3, \dots , a_{n-1}) + 2(-1)^{n-1}
\end{equation}
and
\begin{equation}\label{eqn3}
\owb -a_1, -a_2, \dots , -a_n \cwb = (-1)^n \owb a_1, a_2, \dots ,
a_n \cwb +2(-1)^{n-1}.
\end{equation}

Rewriting equation (\ref{eqn0}) we have 
\begin{equation}\label{eqn5a} 
\Delta_{p/q, \lambda}(x) = \owb a_1, \dots , a_q\cwb + 2( (-1)^q +
   (\lambda/2)^q)
\end{equation}

\begin{notation}\begin{enumerate}
\item
For $0 \leq k \leq n/2$, let $\s{n}{k} = \{ I \subset \{1, 2, \dots
, n\} \mid |I| = n - 2k${\ and $I$ is obtained from $\{ 1, 2, \dots,
n\}$ by deleting $k$ pairs of adjacent elements}$\}$. 
$\s{2k}{k} = \{\emptyset\}$, $\s{2k + 1}{k} = \big\{ \{1\}, \{3\},
\{5\},\dots , \{2k+1\}\big\}$, \dots ,  $\s{n}{0} = \Big\{\{1, 2, 3 , \dots ,
n\}\Big\}$.

\item
For $0 \leq k \leq (n-1)/2$, let $\spr{n}{k} = \{ I \subset \{2, 3, \dots
, n\} \mid |I| = n - 2k -1$ and $I$ is obtained from $\{2, 3, \dots ,
n\}$ by deleting $k$ pairs of adjacent elements$\}$. $\spr{2k + 1}{k} =
\emptyset$, $\spr{2k+2}{k} = \big\{\{2\}, \{4\}, \{6\}, \dots ,
\{2k+2\}\big\}$, \dots ,  $\spr{n}{0} = \{2, 3, \dots , n \}$. 

\item
For $\mathcal{ S}$ a collection of subsets of $\{1, 2, \dots , n - 1\}$
let
$\mathcal{S}\vee \{n\} = \{I \cup \{n\} \mid I \in \mathcal{S} \}$.

\item
For $0 \leq k \leq n/2$, let $\st{n}{k} = \{ I \subset \{1, 2, 3, \dots
, n\} \mid |I| = n - 2k$ and $I$ is obtained from $\{1, 2, \dots , n\}$
by deleting $k$ pairs of adjacent elements, counting $\{n, 1\}$ as an
adjacent pair$\}$. $\st{2k}{k} = \emptyset$, $\st{2k + 1}{k} = \big\{
\{1\}, \{2\}, \dots ,
\{n\}\big\}$, \dots , $\st{n}{0} = \{ 1, 2, 3,
\dots , n\}$. 

\item
For $a_1, a_2, a_3, \dots, a_n$ elements of a commutative ring, and $I =
\{i_1, i_2,\allowbreak i_3, \allowbreak
 \dots , i_k\} \subset \{1, 2, 3, \dots , n\}$, let $a_I
= a_{i_1} a_{i_2} a_{i_3} \cdots a_{i_n}$. We shall adopt the convention
that $a_{\emptyset} = 1$
\end{enumerate}
\end{notation}

Part {(\it ii}) of the next proposition goes back to Sylvester's original
paper on continuants \cite{sylv}; part ({\it iv}) is a straightforward
extension of this. For the reader's convenience we present a proof.

\begin{proposition}\label{prop3}
\begin{enumerate}
\item Suppose $1 \leq k < n/2$, then $\s{n}{k} =\left(\s{n-1}{k} \vee
\{n\} \right) \cup \s{n-2}{k-1}$.

\item 
\[
\sum_{k=0}^{[n/2]} (-1)^k \sum_{I\in \s{n}{k}} a_I 
= a_n \sum_{k=0}^{[(n-1)/2]} (-1)^k \sum_{I\in \s{n-1}{k}} a_I 
-
\sum_{k=0}^{[(n-2)/2]} (-1)^k \sum_{I\in \s{n-2}{k}} a_I.
\]

\item $\st{n}{k} = \s{n}{k} \cup \spr{n-1}{k-1}$ for $1 \leq k \leq n/2$.

\item When $n$ is odd,
\[
\sum_{k=0}^{[n/2]} (-1)^k \sum_{I\in \st{n}{k}} a_I 
= \sum_{k=0}^{[n/2]} (-1)^k \sum_{I\in \s{n}{k}} a_I 
-
\sum_{k=0}^{[(n-1)/2]} (-1)^k \sum_{I\in \spr{n-1}{k}} a_I
\]
When $n$ is even,
\[
\sum_{k=0}^{[n/2]} (-1)^k \sum_{I\in \st{n}{k}} a_I 
= \sum_{k=0}^{[n/2]} (-1)^k \sum_{I\in \s{n}{k}} a_I  
+ (-1)^{n/2} -
\sum_{k=0}^{[(n-1)/2]} (-1)^k \sum_{I\in \spr{n-1}{k}} a_I
\]
\end{enumerate}
\end{proposition}

\proof {\it(i)} Let $I \in \s{n}{k}$. If $n \not\in I$ then $n-1 \not\in
I$ and so $I \in \s{n-2}{k-1}$. Suppose $n \in I$. Let $K= \{1,2,3,
\dots , n\} \setminus I$ and $\dot{I} = I\setminus\{n\}$. Then $\dot{I}
= \{1,2,3, \dots n-1\} \setminus K$; so $\dot{I} \in \s{n-1}{k}$. Hence
$I = \dot{I} \cup \{n\} \in \s{n-1}{k} \vee \{n\}$. 

{\it(ii)} Let us assume that $n =2m$ is even. The same idea works for odd
$n$ but the proof is slightly simpler. Observe
\begin{eqnarray*}
\lefteqn{
\sum_{k=0}^{[n/2]} (-1)^k \sum_{I \in \s{n}{k}} a_I   
=
\sum_{I\in \s{n}{0}} a_I  + \sum_{k=1}^{m-1} (-1)^k \sum_{I \in
\s{n}{k}} a_I  + (-1)^m }\\
&=& \sum_{I\in \s{n}{0}} a_I  + a_n\sum_{k=1}^{m-1} (-1)^k\kern-7pt
\sum_{I \in \s{n-1}{k}} a_I  + \sum_{k=1}^{m-1} (-1)^k\kern-7pt \sum_{I
\in \s{n-2}{k-1}} a_I + (-1)^m \\
&=&
a_n\bigg\{\sum_{I\in \s{n-1}{0}} a_I  + \sum_{k=1}^{m-1}
(-1)^k\kern-7pt \sum_{I \in \s{n-1}{k}} a_I\Bigg\}  + \sum_{k=1}^{m-1}
(-1)^k\kern-7pt \sum_{I
\in \s{n-2}{k-1}} a_I + (-1)^m \\
&=& 
a_n\sum_{k=0}^{m-1} (-1)^k\kern-7pt \sum_{I \in \s{n-1}{k}} a_I
- \sum_{k=0}^{m-1} (-1)^k\kern-7pt \sum_{I
\in \s{n-2}{k}} a_I \\
&=& 
a_n\sum_{k=0}^{[(n-1)/2]} (-1)^k\kern-7pt \sum_{I \in \s{n-1}{k}} a_I
- \sum_{k=0}^{[(n-2)/2]} (-1)^k\kern-7pt \sum_{I
\in \s{n-2}{k}} a_I.
\end{eqnarray*}

{\it (iii)} For $I \in \st{n}{k}$ let $K_1= \{1,2,3, \dots , n\}
\setminus I$ and $K_2= \{2,3, \dots , n\}
\setminus I$. $\min\{i \mid i \in K_1\}$ is odd if and only if $I \in
\s{n}{k}$ and  $\min\{i \mid i \in K_2\}$ is even if and only if $I \in
\spr{n-1}{k-1}$

{\it (iv)} Suppose $n =2m$.  Then
\begin{eqnarray*}
\lefteqn{
\sum_{k=0}^{[n/2]} (-1)^k \sum_{I\in \st{n}{k}} a_I 
= \sum_{I\in \st{n}{0}} a_I +
\sum_{k=1}^{m-1} (-1)^k \sum_{I\in \st{n}{k}} a_I 
+ (-1)^{m} \sum_{I\in \st{n}{m}} a_I }\\ 
&=& \Bigg(\sum_{I\in \s{n}{0}} a_I +
\sum_{k=1}^{m} (-1)^k\kern-7pt \sum_{I\in \s{n}{k}} a_I\Bigg) +
\sum_{k=1}^{m} (-1)^k\kern-7pt \sum_{I\in \spr{n-1}{k-1}} a_I + (-1)^m\\
&=&
\sum_{k=0}^{m} (-1)^k \sum_{I\in \s{n}{k}} a_I 
- \sum_{k=0}^{m-1} (-1)^k\kern-7pt \sum_{I\in \spr{n-1}{k}} a_I + (-1)^m.
\end{eqnarray*}
The case of $n$ odd is similar. \qed

\begin{corollary}\label{cor4}
Let $a_1, a_2, \dots , a_n$ be elements of a commutative ring.
\begin{enumerate}
\item $\displaystyle{
(a_1, \dots , a_n) = \sum_{k=0}^{[n/2]} (-1)^k \sum_{I \in \s{n}{k}} a_I
}$
\item $$
\owb a_1, \dots , a_n\cwb 
= \begin{cases}
\displaystyle{
\sum_{k=0}^{[n/2]} (-1)^k \sum_{I \in \st{n}{k}} a_I} + 2
& n \textrm{ odd} \\
\displaystyle{
\sum_{k=0}^{[n/2]} (-1)^k \sum_{I \in \st{n}{k}} a_I - 2 +
(-1)^{n/2}} 
& n \textrm{ even} \\ \end{cases}$$
\end{enumerate}
\end{corollary}

\proof {\it(i)} For $n =1$ the left hand side and the right hand side
equal $a_1$. Both sides satisfy the same recurrence relation. 

{\it (ii)}  By equation (\ref{eqn5a}) 
$$
\owb a_1, \dots , a_n \cwb = (a_1, \dots , a_n) - (a_2, \dots , a_{n-1})
- (-1)^n 2
$$
so the result now follows from ({\it i}) and Proposition
\ref{prop3} ({\it iii}).
\qed

\begin{proposition}\label{prop5}
Let $1 \leq p < q$ be relatively prime, $\theta = p/q$, and $a_k = 2
\cos(2 \pi k \theta)$. Then
$$
a_1 a_2 \cdots a_q =
\begin{cases}
\phantom{-}0 & q \equiv 0 \pmod 4 \\
\phantom{-}2 & q \equiv 1, 3 \pmod 4 \\
-4 & q \equiv 2 \pmod 4. \\
\end{cases}$$
\end{proposition}

\proof Let $T_q$ be the $q$th Chebyshev polynomial of the first kind.
The constant term of $T_q$ is 0 for $q$ odd and $(-1)^{q/2}$ for $q$
even. The result now follows from the identity (see e.g.
\cite[\S 1.2]{rivlin})
$$ \prod_{i=1}^q (x - a_i) = 2( T_q(x/2) -1).$$\qed

The statement of the main theorem follows. Its proof will be given at the end of 
the next section.

\begin{theorem}\label{main}
Let $1 \leq p < q$ be relatively prime, $a_k = 
\allowbreak 2  \cos(2 \pi k \theta)$, and $\theta = p/q$. For $1
\leq k < q/2$, 
$$
\sum_{I\in \st{q}{k}} a_I = 0.$$
\end{theorem}

\begin{corollary}
Let $1 \leq p < q$ be relatively prime, $\theta = p/q$, $\lambda > 0$,
and $a_k = \lambda \cos(2 \pi k \theta)$. Then
$$
\owb a_1, a_2, \dots , a_q\cwb
= \begin{cases}
\phantom{-}0 & q \equiv 0 \pmod 4 \\
2(1 + (\lambda/2)^q) & q \equiv 1, 3 \pmod 4 \\
-4(1 + (\lambda/2)^q) & q \equiv 2 \pmod 4 \\ \end{cases}
$$
and $\Delta_{\theta, \lambda}(0) = (-1)^{q/2} 2 ( 1 + (\lambda/2)^q)$ for
$q$ even.
\end{corollary}

\proof 
Suppose $q$ is even. By Theorem \ref{main} all the terms of
$$
\sum_{k=0}^{[q/2]}(-1)^k \sum_{I\in \st{q}{k}} a_I
$$
are zero except the terms for $k = 0$ and $k= q/2$. The term for $k = 0$
is $a_1 a_2 \cdots a_q$. The term for $k = q/2$ is $(-1)^{q/2}$. Thus
when $q= 4m$ we have by Proposition \ref{prop5}
$$
\owb a_1 , a_2 , \dots , a_q\cwb = a_1 a_2 \cdots a_q - (-1)^q 2 +
(-1)^{q/2}2 = 0,
$$
and when $q = 4m +2$,
\begin{eqnarray*}
\owb a_1 , a_2 ,\dots , a_q\cwb = 
a_1 a_2 \cdots a_q - (-1)^q 2 + (-1)^{q/2}2 =
-4(1 + (\lambda/2)^q).
\end{eqnarray*}
To obtain the final claim we apply equation (\ref{eqn5a}).\qed

From the corollary and equation (\ref{eqn00}) we have the theorem which
corrects an error in \cite[p. 232]{choi-elliott-yui}.

\begin{theorem}
The coefficient of $x^0$ in $\det(xI_q - h_{\theta, \lambda}(z_1, z_2))$ is 
$$ \mbox{} - (z_1^q + z_1^{-q} + (\lambda/2)^q ( z_2^q + z_2^{-q})) +
(-1)^{q/2} 2 (1 + (\lambda/2)^q))$$ when $q$ is even and  
$ - (z_1^q + z_1^{-q} + (\lambda/2)^q ( z_2^q +
z_2^{-q}))$ when $q$ is odd.
\end{theorem}

\section{Proof of the Main Theorem}

\begin{theorem}\label{thm9}
Suppose $a_1, a_2, \dots , a_q$ are elements in a commutative
ring and let $a_{q+1} = a_1$. For $I \subset \{ 1, 2, \dots , q\}$, let
$I^c = \{ 1, 2, \dots , q\} \setminus I$, be the complement of $I$ in $\{
1, 2, \dots , q\}$. Then
\begin{eqnarray*}\lefteqn{
\sum_{I\in \st{q}{k}} a_{I^c} =
\sum_{i_1 = 1}^{q - 2(k-1)} 
\sum_{i_2 = i_1 + 2}^{q - 2(k-2)} 
\cdots
\sum_{i_k = i_{k-1}+2}^q \ 
\prod_{j=1}^k a_{i_j}a_{i_j+1}} \\
&&\mbox{} -
a_1a_2\Bigg[
\sum_{i_1 = 3}^{q - 2(k-2)} 
\sum_{i_2 = i_1 + 2}^{q - 2(k-3)} 
\cdots  
\sum_{i_{k-2} = i_{k-3}+2}^{q-2} \ 
\prod_{j=1}^{k-2} a_{i_j}a_{i_j+1} \Bigg]
a_q a_{q+1}.
\end{eqnarray*}
\end{theorem}

\proof Recall that elements of $\st{q}{k}$ are obtained by deleting $k$
adjacent pairs $\{i, i +1\}$ from $\{ 1, 2, \dots , q\}$, counting $q$
and 1 as adjacent.  So if $I^c \in \st{q}{k}$ then $I = \{i_1, j_1,
i_2, j_2, \dots , i_k , j_k\}$ with $1 \leq i_1, j_1 = i_1 + 1 < i_2,
\dots , j_{k-1} = i_{k-1} + 1 <  i_k \leq  q$ and either $j_k = i_k + 1$
if $i_k < q$ or $j_k = 1$ if $i_k = q$. 

Now let $\ts{q}{k} = \{\, \{i_1, j_1, i_2, j_2, \dots , i_k, j_k\} \mid
1 \leq i_1, j_1 = i_1 + 1 < i_2, \dots ,\allowbreak
 j_{k-1} = i_{k-1} + 1 < i_k
\leq q, j_k = i_k +1 \, \}$. Define $\phi : \{ 1, 2, \dots , q, q + 1 \}
\rightarrow \{1, 2, \dots , q\}$ by $\phi(q+1) = 1$ and $\phi(i) = i$
for $i \leq q$. Then $a_{\phi(I)} = a_I$ for $I \in \ts{q}{k}$. 

If $I = \{i_1, j_1, i_2, j_2, \dots , i_k, j_k\}$ and $i_1 =1$ and $i_k =
q$ then $\phi(I)^c \not \in \st{q}{k}$ because $\phi(j_k) = \phi(i_1) =
1$ and the pairs must be disjoint. So let $\tpr{q}{k} = \{ \, \{1, 2,
i_1, j_1, \dots ,\ab i_{k-1}, j_{k-1}, q , q+1\} \mid 3 \leq i_1, j_1 =
i_1 + 1 < i_2, \dots , i_{k-1} \leq q-2, j_{k-1} = i_{k-1} + 1\, \}$.

For $I \in \ts{q}{k} \setminus \tpr{q}{k}$, $\phi(I)^c \in \st{q}{k}$
and $\phi: \ts{q}{k} \setminus \tpr{q}{k} \rightarrow \st{q}{k}$ is a
bijection. This with the identity $a_{\phi(I)} = a_I$ proves the
theorem. \qed

\begin{lemma}\label{lemma9}
\begin{enumerate}
\item    For $q \geq 1$,
$$
\sum_{i=1}^q (x^{-i} + x^i)^{-1}
             (x^{-i-1} + x^{i+1})^{-1} 
= \frac{ x^{-q} - x^q }%
{(x^{-2} - x^2 ) 
                       (x^{-q-1} + x^{q+1}) }. $$
\item For $k \geq 1$
$$\prod_{i=1}^{2k} (x^{-i} + x^i)^{-1} =
 \prod_{i=1}^k \frac{  \displaystyle (x^{-i} - x^{i}) }
{(x^{-2i} - x^{2i}) 
(x^{-(k+i)} + x^{k+i}) }. $$
\end{enumerate}
\end{lemma}

\proof {\it (i)} One checks directly that the
formula holds when $q=1$, then {\it (i)} follows by induction on $q$.

{\it (ii)} follows from the identity
$$ 
\frac{x^{-i} - x^i }{ (x^{-2i} - x^{2i}) (x^{-k-i} + x^{k+i})}
=
\frac{1 }{
(x^{-i} + x^{i}) (x^{-k-i} + x^{k+i}) 
}$$ \qed

\begin{corollary}\label{cor11}
For $q \geq 5$
$$\frac{(x^{-1} + x)^{-1} (x^{-2} + x^2)^{-1}}
{(x^{-q} + x^q) (x^{-q-1} + x^{q+1})}
\sum_{i=3}^{q-2}
(x^{-i} + x^{i})^{-1} (x^{-i-1} + x^{i+1})^{-1} $$
$$ = \frac{
(x^{-3} - x^3) (x^{-q+4} - x^{q-4})  }{
(x^{-4} - x^4) (x^{-6} - x^6) (x^{-q+1} + x^{q-1}) (x^{-q} + x^q) (x^{-q-1} +
x^{q+1})  }$$
\end{corollary}

\proof By Lemma \ref{lemma9} {\it (i)} 
\begin{eqnarray*}\lefteqn{
\sum_{i=3}^{q-2} (x+{-i} + x^i)^{-1} (x^{-i-1} + x^{i+1})^{-1} } \\
&=& \frac{
x^{-q+2} -x ^{q-2} }{  (x^{-2} - x^2) (x^{-q+1} + x^{q-1}) } - 
\frac{x^{-2} - x^2  }{  (x^{-2} - x^2) (x^{-3} + x^{3})} \\
&=& \frac{(x^{-q+4} - x^{q-4}) (x^{-1} + x)  }{
(x^{-2} - x^2) (x^{-3} + x^3) (x^{-q+1} + x^{q-1})}
\end{eqnarray*}
The result then follows by multiplying both sides by $(x^{-1} + x) (x^{-2} +
x^2) (x^{-3} + x^3) (x^{-q+1} + x^{q-1}) $. \qed

\begin{theorem} \label{thm12}
For $k \geq 1$ and $q \geq 2k-1$,
$$
\sum_{i_1 = 1}^{q - 2(k-1)} 
\sum_{i_2 = i_1 + 2}^{q - 2(k-2)} 
\cdots
\sum_{i_k = i_{k-1}+2}^q \ 
\prod_{j=1}^k (x^{-i_j} + x^{i_j})^{-1}
              (x^{-i_j-1} + x^{i_j+1})^{-1}
$$
\begin{equation}\label{eqn6}
\mbox{} = \frac{
\prod_{i=k-1}^{2k-2} (x ^{-(q-i)} - x^{q-i})
}{ 
\prod_{i=1}^k (x^{-2i} - x^{2i})
\prod_{i=-1}^{k-2} (x^{-(q-i)} + x^{q-i})
}\end{equation}

\end{theorem}

\proof We prove the equation by induction on $k$. When
$k=1$ the equation holds by Lemma \ref{lemma9} {\it (i)}. 
Lemma \ref{lemma9} {\it (ii)}
shows that for arbitrary $k$ the formula holds for $q=
2k-1$; so we fix $k$ and proceed by induction on $q$.
Let $S_{k,q}$ and $T_{k,q}$ denote respectively the left hand  and right
hand sides of equation (\ref{eqn6}). 

If we write $S_{k,q}$ as a sum of two terms, the first in which $i_k < q$ and
the second when $i_k = q$, we see that $S_{k,q}$ satisfies the recurrence
relation
$$ S_{k,q} = S_{k, q-1} + (x^{-q} + x^q)^{-1}
(x^{-q-1} + x^{q+1})^{-1} S_{k-1, q-2}
$$

Thus we have only to show that $T_{k,q}$ satisfies the same relation. Now
$$
T_{k, q-1} = 
\frac{
\prod_{i=k}^{2k-1} (x^{-(q-i)} - x^{q-i}) }{
\prod_{i=1}^k (x^{-2i} - x^{2i})
\prod_{i=0}^{k-1} (x^{-(q-i)} + x^{q-i})
}$$
and
$$
T_{k-1, q-2} = 
\frac{
\prod_{i=k}^{2k-2} (x^{-(q-i)} - x^{q-i}) }{
\prod_{i=1}^{k-1} (x^{-2i} - x^{2i})
\prod_{i=1}^{k-1} (x^{-(q-i)} + x^{q-i})
}$$
The proof of the recurrence relation for $T_{k,q}$ is thus reduced to verifying
that
$$ \frac{
(x^{-(q-(k-1)} - x^{q - (k-1)}) (x^{-(q - (k-1)} + x^{q - (k-1)})
}{
(x^{-q-1} + x^{q+1}) (x^{-q} + x^q)}
$$
$$ \mbox{} =
\frac{x^{-(q-(2k-1))} - x^{q - (2k-1)} }{
x^{-q} + x^q}
+
\frac{x^{-2k} - x^{2k} }{
(x^{-q} + x^q) (x^{-q-1} + x^{q+1})}
$$
\qed

\begin{theorem}\label{thm13}
For $k \geq 3$ and $q \geq 2k-1$,
{%
\begin{eqnarray}
\label{eqn7}\lefteqn{
(x^{-1} + x^1)^{-1} (x^{-2} + x^2)^{-1} 
(x^{-q} + x^q)^{-1} (x^{-q-1} + x^{q+1})^{-1}  } \nonumber \\ 
& &
\times \sum_{i_1 = 3}^{q - 2(k-2)} 
\sum_{i_2 = i_1 + 2}^{q - 2(k-3)} 
\cdots  
\sum_{i_{k-2} = i_{k-3}+2}^{q-2} \ 
\prod_{j=1}^{k-2} (x^{-i_j} + x^{i_j})^{-1}
              (x^{-i_j-1} + x^{i_j+1})^{-1}  \nonumber \\
& = &
\frac{(x^{-k+1} - x^{k-1}) (x^{-k} - x^k) 
\prod_{i=k+1}^{2k-2}
(x^{-q+i} - x^{q -i} ) 
}{ 
\prod_{i=1}^k
(x^{-2i} - x^{2i})
\prod_{i=-1}^{k-2}(x^{-q+i} + x^{q - i} ) } 
\end{eqnarray}}
\end{theorem}

\proof Let us denote the left and right hand sides of the identity by
$S_{k,q}$ and $T_{k,q}$ respectively. By Corollary \ref{cor11} $S_{3,q} =
T_{3,q}$. We write $S_{k,q}$ as the sum of two terms: in the first $i_{k-2} <
q - 2$ and in the second $i_{k-2} - q-2$. As in the proof of the previous
theorem we obtain a recurrence relation, in this case:
\begin{eqnarray*}\lefteqn{
S_{k,q} = S_{k, q-1} (x^{-q-1} + x^{q+1})^{-1} (x^{-q+1} + x^{q-1}) }\\
&&\mbox{} +
S_{k-1, q-2} (x^{-q} + x^q)^{-1} (x^{-q-1} + x^{q+1})^{-1}
\end{eqnarray*}
It is routine to verify that $T_{k,q}$ satisfies the same recurrence
relation. \qed

Subtracting equation  (\ref{eqn7})  from equation (\ref{eqn6}) yields.

\begin{corollary} \label{cor14}
\begin{eqnarray}\label{eqn8}
\lefteqn{
\sum_{i_1 = 1}^{q - 2(k-1)} 
\sum_{i_2 = i_1 + 2}^{q - 2(k-2)} 
\cdots
\sum_{i_k = i_{k-1}+2}^q \ 
\prod_{j=1}^k (x^{-i_j} + x^{i_j})^{-1}
              (x^{-i_j-1} + x^{i_j+1})^{-1}}\nonumber \\ 
&&\mbox{}-
(x^{-1} + x^1)^{-1} (x^{-2} +
x^2)^{-1} (x^{-q} + x^q)^{-1} (x^{-q-1} + x^{q+1})^{-1} \nonumber \\
&\times &
\sum_{i_1 = 3}^{q - 2(k-2)} 
\sum_{i_2 = i_1 + 2}^{q - 2(k-3)} 
\cdots\sum_{i_{k-2} = i_{k-3}+2}^{q-2} \ 
\prod_{j=1}^{k-2} (x^{-i_j} + x^{i_j})^{-1}
              (x^{-i_j-1} + x^{i_j+1})^{-1} \nonumber\\
&=&
\frac{ (x^q - x^{-q}) \displaystyle
\prod_{i=k+1}^{2k-2} (x^{-q+i} - x^{q-i}) }{
\displaystyle \prod_{k=1}^q (x^{2i} - x^{-2i}) \quad
\prod_{i=-1}^{k-2}(x^{-q+i} + x^{q-i}) }
\end{eqnarray}

\end{corollary}

\bigskip

\noindent{\it \label{mainproof}
Proof of theorem \ref{main} }  
We recall that $1 \leq p < q$ and $p$ and $q$ are relatively prime.
We set $\theta = p/q$ and $a_j = 2 \cos(2 \pi j \theta)$. We shall split
the proof into two cases.

\medskip\noindent{\it
Case 1: $q \not\equiv 0 \pmod 4$}. When $q\not\equiv 0 \pmod 4$  $a_j
\not = 0$ for all $j$; moreover when $x = e^{2 \pi i \theta}$, $x^{4i}
\not = 1$ and $x^{2(q - i)}
\not = -1$ for all $i$. Thus the denominator on the right hand side of
(\ref{eqn8}) does not vanish but the numerator does. Hence by Theorem
\ref{thm9}
\[ \sum_{I \in \st{q}{k}} (a_{I^c})^{-1} = 0 \]
Upon multiplying by $a_1 a_2 \cdots a_q$ we obtain that

\[
\sum_{I \in \st{q}{k}} a_I = a_1a_2\cdots a_q 
\sum_{I \in \st{q}{k}} (a_{I^c})^{-1} = 0.
\]

\medskip\noindent{\it
Case 2: $q \equiv 0 \pmod 4$}. 
Again we wish to show that $\sum_{I \in  \st{q}{k}} a_{I} = 0$ and so  we
must multiply both sides of equation (\ref{eqn8}) by $\prod_{i=1}^q
(x^{-i} + x^i)$ and evaluate at $x = e^{2 \pi i \theta}$. 

The denominator of the right hand side of
(\ref{eqn8}) is zero when $x^{4q} = 1$ or $x^{2(q-j)} = -1$, i.e. when
$i= j = q/4$; the corresponding factors are $x^{-q/2} - x^{q/2}$ and
$x^{-3q/4} + x^{3q/4}$ respectively. 

Apart from the factor $x^{-q} - x^q$, the numerator of the right hand
side of equation (\ref{eqn8}) is zero only when $x^{2(q-i)} = 1$, i.e.
when $i = q/2$. This produces the factor $x^{-q/2} - x^{q/2}$ which
cancels one of the zeros in the denominator. The other zero is
cancelled when we  multiply by  $\prod_{i=1}^q (x^{-i} + x^i)$. Hence
the product of $\prod_{i=1}^q (x^{-i} + x^i)$ and the right side of
(\ref{eqn8}) is zero when $x = e^{2 \pi i \theta}$. \qed


\end{document}